\documentclass[11pt]{amsart}         
\DeclareMathOperator{\Hom}{Hom}
\DeclareMathOperator{\Com}{Com}
\DeclareMathOperator{\Cup}{Cup}
\DeclareMathOperator{\I}{I}
\DeclareMathOperator{\1}{id}

\DeclareMathOperator{\dev}{dev}
\DeclareMathOperator{\x}{\times}
\newcommand{\NN}{\mathbb{N}}
\newcommand{\EEnd}{\mathcal End}
\newcommand{\EE}{\mathcal E}
\newcommand\bul{\bullet}
\newcommand{\de}{\delta}
\newcommand{\Ga}{\Gamma}
\renewcommand{\=}{:=}
\renewcommand{\t}{\otimes}
\renewcommand{\o}{\circ}

\renewcommand{\u}{\smile}
\newtheorem{thm}{Theorem}[section]
\newtheorem{prop}[thm]{Proposition}

\newtheorem{lem}[thm]{Lemma}
\newtheorem*{Boundary Lemma}{Boundary Lemma}
\newtheorem*{Main Theorem}{Main Theorem}
\theoremstyle{definition}
 \newtheorem{defn}[thm]{Definition}
\theoremstyle{remark}
 \newtheorem{rem}[thm]{Remark}
\theoremstyle{remark}
 \newtheorem{rems}[thm]{Remarks}
\theoremstyle{definition}
\newtheorem{exam}[thm]{Example}
\theoremstyle{definition}
 \newtheorem{recap}[thm]{Recapitulation}
\theoremstyle{definition}
 \newtheorem{notations}[thm]{Notations}
\numberwithin{equation}{subsection}

\begin{document}
\baselineskip0.95\baselineskip
\null
\vskip2.5true cm

\begin{center}
{\large\bf TETRACOMPOSITION}\\
\vskip22pt
{\large Liivi Kluge and Eugen Paal}
\end{center}

\begin{quote}
\item
{\bf Abstract.}
We consider basic algebraic constructions associated with an \emph{abstract}
pre-operad, such as a $\smile$-algebra, total composition $\bul$,
\emph{pre-cobound\-ary} operator $\de$, \emph{tribraces} $\{\cdot,\cdot,\cdot\}$
and \emph{tetrabraces} $\{\cdot,\cdot,\cdot,\cdot\}$.
A derivation deviation of the
pre-coboundary operator over the tetrabraces is calculated
in terms of the $\smile$-multiplication and tribraces.

\smallskip
{\bf Classification~(MSC2000).} 18D50.

\end{quote}

\renewcommand{\thefootnote}{}
\footnote{math.QA/0105231}
\footnote{Research supported in part by the ESF grant 3654.}

\section{Introduction and outline of the paper}

\noindent We consider basic algebraic constructions associated
with an \emph{abstract} pre-operad $C$ (Sec.~\ref{pre-operad}),
such as a $\smile$-algebra (Sec.~\ref{cup-mul}), total
composition $\bul$ (Sec.~\ref{cup-mul}), \emph{pre-coboundary}
operator $\de$ (Sec.~\ref{cup and pre-coboundary}),
\emph{tribraces} $\{\cdot,\cdot,\cdot\}$ (Sec.~\ref{tribraces})
and \emph{tetrabraces} $\{\cdot,\cdot,\cdot,\cdot\}$
(Sec.~\ref{tetrabraces}). Main result of the present paper is the
Main Theorem of Sec.~\ref{tetrabraces}. By defining (see degree
notations in Sec.~\ref{pre-operad}) a \emph{derivation deviation}
(Sec.~\ref{tetrabraces}) of the pre-coboundary operator $\de$ over the
tetrabraces,
\begin{align*}
\dev_{\{\cdot,\cdot,\cdot,\cdot\}}\de\,(h\t f\t g\t b)
    &\=\de\{h,f,g,b\}-\{h,f,g,\de b\}-(-1)^{|b|}\{h,f,\de g,b\}\\
    &\quad-(-1)^{|b|+|g|}\{h,\de f,g,b\}
        -(-1)^{|b|+|g|+|f|}\{\de h,f,g,b\},
\end{align*}
the Main Theorem states that in a the pre-operad $C$ one has
\begin{align*}
(-1)^{|b|}\dev_{\{\cdot,\cdot,\cdot,\cdot\}}\de\,(h\t f\t g\t b)
    &=\{h,f,g\}\u b-\{h,f,g\u b\}\\
    &-(-1)^{|g|}\{h,f\u g,b\}+(-1)^{|h|f+|g|}f\u\{h,g,b\}.
\end{align*}
Our result is similar to that of \cite{GerVor94,VorGer}, where it was
announced that for the \emph{Hochschild cochains} the extra terms with $\de
h$, $\de f$, $\de g$ and $\de b$  also appear. We work out these extra terms
for an \emph{abstract} pre-operad and give their interpretation via a
\emph{derivation deviation} of $\de$ over the tetrabraces.

The paper can be seen as an an extension of \cite{KPS} and \cite{KP}.

\section{Pre-operad (composition system)}
\label{pre-operad}

\noindent Let $K$ be a unital commutative associative ring, and let $C^n$
($n\in\NN$) be unital $K$-modules. For \emph{homogeneous} $f\in C^n$, we
refer to $n$ as the \emph{degree} of $f$ and often write (when it does not
cause confusion) $f$ instead of $\deg f$; for example $(-1)^f\=(-1)^n$,
$C^{f}\=C^{n}$ and $\o_f\=\o_n$. Also, it is convenient to use the shifted
(\emph{desuspended}) degree $|f|\=n-1$. Throughout this paper, we assume
that $\t\=\t_K$ and work with \emph{homogeneous} elements of $C$.

\begin{defn}
A linear (right) \emph{pre-operad} (\emph{composition system})
with coefficients in $K$ is a sequence $C\=\{C^n\}_{n\in\NN}$ of unital
$K$-modules (an $\NN$-graded $K$-module),
such that the following conditions hold:
\begin{enumerate}
\item
For $0\leq i\leq m-1$ there exist \emph{partial compositions}
\[
  \o_i\in\Hom\,(C^m\t C^n,C^{m+n-1}),\qquad |\o_i|=0.
\]
\item
For all $h\t f\t g\in C^h\t C^f\t C^g$, the \emph{composition relations} hold,
\[
(h\o_i f)\o_j g=
   \begin{cases}
    (-1)^{|f||g|} (h\o_j g)\o_{i+|g|}f
                       &\text{if $0\leq j\leq i-1$},\\
    h\o_i(f\o_{j-i}g)  &\text{if $i\leq j\leq i+|f|$},\\
    (-1)^{|f||g|}(h\o_{j-|f|}g)\o_i f
                       &\text{if $i+f\leq j\leq|h|+|f|$}.
\end{cases}
\]
\item
There exists a unit $\I\in C^1$ such that
\[
\I\o_0 f=f=f\o_i \I,\qquad 0\leq i\leq |f|.
\]
\end{enumerate}

In the 2nd item, the \emph{first} and \emph{third} parts of the
defining relations turn out to be equivalent.
\end{defn}

\begin{rems}
A pre-operad is also called a
              \emph{comp(osition) algebra}
           or \emph{asymmetric operad}
           or \emph{non-symmetric operad}
           or \emph{non-$\Sigma$ operad}.
The concept of (\emph{symmetric}) \emph{operad} was formalized by
P. May \cite{may72} as a tool for the theory of iterated loop spaces.
Recent studies and applications can be found in \cite{Rene}.

Above we modified the Gerstenhaber \emph{comp algebra} defining relations
\cite{CGS93,GGS92Am} by introducing the sign $(-1)^{|f||g|}$ in the
defining relations of the pre-operad.
The modification enables us to keep track of (control) sign changes more
effectively.
One should also note that (up to sign) our $\o_i$ is Gerstenhaber's $\o_{i+1}$
from \cite{GGS92Am,CGS93}; we use the original (non-shifted)
convention from \cite{Ger,Ger68}.
\end{rems}

\begin{exam}[endomorphism pre-operad]
\label{HG}

Let $A$ be a unital $K$-module and
$\EE_A^n\={\EEnd}_A^n\=\Hom\,(A^{\t n},A)$.
Define (cf. \cite{Ger,Ger68,GGS92Am}) the partial compositions
for $f\t g\in\EE_A^f\t\EE_A^g$ as
\[
f\o_i g\=(-1)^{i|g|}f\o(\1_A^{\t i}\t g\t\1_A^{\t(|f|-i)}),
         \qquad 0\leq i\leq |f|.
\]
Then $\EE_A\=\{\EE_A^n\}_{n\in\NN}$ is a pre-operad
(with the unit $\1_A\in\EE_A^1$) called the \emph{endomorphism pre-operad}
of $A$.
A few examples (without the sign factor) can be found in
\cite{Ger68,GGS92Am} as well.
We use the original indexing of \cite{Ger,Ger68} for the
defining formulae.
\end{exam}

\begin{exam}[associahedra]
A geometrical  example of a pre-operad is provided by the Stasheff
\emph{associahedra}, which was first constructed in \cite{Sta63}. Quite a
surprising realization of the associahedra as \emph{truncated simplices}
was discovered and studied in \cite{ShnSte94,Sta97,Markl97}.
\end{exam}

\begin{notations}[scope of a pre-operad]
\label{scope}
The {\it scope} of $(h\o_i f)\o_j g$ is given by
\[
0\leq i\leq|h|,\qquad 0\leq j\leq|f|+|h|.
\]
It follows from the defining relations of a pre-operad that
the scope is a disjoint union of
\allowdisplaybreaks
\begin{align*}
B &\=\{(i,j)\in \NN\x\NN\,|\, 1\leq i\leq|h|  \,;\, 0   \leq j\leq i-1\},\\
A &\=\{(i,j)\in \NN\x\NN\,|\, 0\leq i\leq|h|  \,;\, i   \leq j\leq i+|f|\},\\
G &\=\{(i,j)\in \NN\x\NN\,|\, 0\leq i\leq|h|-1\,;\, i+f \leq j\leq |f|+|h|\}.
\end{align*}
\allowdisplaybreaks
Note that the triangles $B$ and $G$ are symmetrically situated
with respect to the parallelogram $A$ in the scope $B\sqcup A\sqcup G$.
The (recommended and impressive) picture is left for a reader as
an exercise.
\end{notations}

\begin{recap}
The defining relations of a pre-operad read
\[
(h\o_i f)\o_j g=
\begin{cases}
   (-1)^{|f||g|}(h\o_j g)\o_{i+|g|}f &\text{if $(i,j)\in B$},\\
   h\o_i(f\o_{j-i}g)                 &\text{if $(i,j)\in A$},\\
   (-1)^{|f||g|}(h\o_{j-|f|}g)\o_i f &\text{if $(i,j)\in G$}.
\end{cases}
\]
The \emph{first} ($B$) and \emph{third} ($G$) parts of the relations
are equivalent.
\end{recap}

\section{Cup and total composition}
\label{cup-mul}

\noindent
In Sections~\ref{cup-mul}-\ref{tribraces}, we recall from \cite{KPS,KP}
main facts about cup-multiplication, total composition, pre-coboundary operator
and tribraces in an abstract pre-operad.
Fix $\mu\in C^2$.

\begin{defn}
In a pre-operad $C$, the \emph{cup-multiplication}
$\u\=\u_\mu\:C^f\t C^g\to C^{f+g}$ is defined
(cf. \cite{Ger,Ger68,GGS92Am}) by
\[
f\u g\=(-1)^f(\mu\o_0 f)\o_f g\in C^{f+g},
\qquad|\smile|=1,\qquad f\t g\in C^f\t C^g.
\]
Note that $\o_f\=\o_{\deg f}$, and $|\smile|=1$ means that $\deg\smile=0$.
The pair $\Cup C\=\{C,\u\}$ is called a $\u$-algebra of $C$.
\end{defn}

\begin{exam}
For the endomorphism pre-operad (Example \ref{HG}) $\EE_A$,
one has
\[
f\u g=(-1)^{fg}\mu\o(f\t g)\in\EE_A^{f+g},
      \qquad \mu\t f\t g\in \EE_A^2\t \EE_A^f\t \EE_A^g.
\]
\end{exam}

\begin{prop}
\label{cuppro}
In a pre-operad $C$, one has
\[
\mu\o_0 f=(-1)^f f\u\I,
\quad \mu\o_1 f=-\I\u f,
\quad f\u g=-(-1)^{|f|g}(\mu\o_1 g)\o_0f.
\]
\end{prop}

\begin{lem}[{\rm\cite{KPS,KP}}]
\label{cup}
In a pre-operad $C$, the following composition relations hold:
\[
(f\u g)\o_j h=
   \begin{cases}
     (-1)^{g|h|}(f\o_j h)\u g &\text{if\,\,\, $0\leq j\leq|f|$},\\
     f\u(g\o_{j-f}h)          &\text{if\,\,\,  $f\leq j\leq|g|+f$}.
   \end{cases}
\]
\end{lem}

\begin{defn}
In a pre-operad $C$, the \emph{total composition}
$\bul\:C^f\t C^g\to C^{f+g-1}$ is defined (cf. \cite{GGS92Am,CGS93}) by
\[
f\bul g\=\sum_{i=0}^{|f|}f\o_i g\in C^{f+g-1},
     \qquad |\bul|=0, \qquad f\t g\in C^f\t C^g.
\]
The pair $\Com C\=\{C,\bul\}$ is called the \emph{composition algebra} of $C$.
\end{defn}

\begin{thm}[{\rm\cite{KPS,KP}}]
\label{right der}
In a pre-operad $C$, one has
\[
(f\u g)\bul h=f\u(g\bul h)+(-1)^{|h|g}(f\bul h)\u g.
\]
\end{thm}

\begin{rem}
This theorem tells us that \emph{right} translations in $\Com C$ are
(right) derivations of the $\u$-algebra.
In turns out that the \emph{left}
translations in $\Com C$ are not derivations of the $\u$-algebra
(see Theorem \ref{tri-main}).
\end{rem}

\section{Pre-coboundary operator}
\label{cup and pre-coboundary}

\begin{defn}
In a pre-operad $C$, define a \emph{pre-coboundary} operator
$\de\=\de_\mu$ by
\[
-\de f\=[f,\mu]\=f\bul\mu-(-1)^{|f|}\mu\bul f,
             \qquad \mu\t f\in C^2\t C^f.
\]
\end{defn}

\begin{exam}
In the Gerstenhaber theory \cite{Ger},
$\de_\mu$ is the Hochschild \emph{coboundary operator} with
the property $\de^2_\mu=0$, the latter is due to the \emph{associativity}
$\mu\bul\mu=0$.

In this paper, we do not assume the associativity constraint for $\mu$,
thus in general $\de^2_\mu\neq0$.
\end{exam}

\begin{prop}
In a pre-operad $C$, one has
\[
-\de f=f\u\I+f\bul\mu+(-1)^{|f|}\I\u f.
\]
\end{prop}

\begin{defn}
The \emph{derivation deviation} of $\de$ over $\bul$ is defined by
\[
\dev_\bul\de\,(f\t g)
   \=\de(f\bul g)-f\bul\de g-(-1)^{|g|}\de f\bul g.
\]
\end{defn}

\begin{thm}[{\rm\cite{KPS}}]
\label{bul-main}
In a pre-operad $C$, one has
\[
(-1)^{|g|}\dev_\bul\de\,(f\t g)=f\u g-(-1)^{fg}g\u f.
\]
\end{thm}

\section{Tribraces}
\label{tribraces}

\begin{defn}[Gerstenhaber tribraces {\rm(cf. \cite{Ger,CGS93})}]
Let $h\t f\t g\in C^h\t C^f\t C^g$.
The Gerstenhaber \emph{tribraces} $\{\cdot,\cdot,\cdot\}$ are
defined as a double sum over the triangle $G$ by
\[
\{h,f,g\}\=\sum_{(i,j)\in G}(h\o_i f)\o_j g\in C^{h+f+g-2},
    \qquad |\{\cdot,\cdot,\cdot\}|=0.
\]
\end{defn}

\begin{thm}[Getzler identity {\rm\cite{KPS}}]
\label{Getzler}
In a pre-operad $C$, one has
\begin{align*}
(h,f,g)&\=(h\bul f)\bul g-h\bul(f\bul g)\\
       &\,\,=\{h,f,g\}+(-1)^{|f||g|}\{h,g,f\}.
\end{align*}
\end{thm}

\begin{thm}[Gerstenhaber identity]
\label{Gerst}
In a pre-operad $C$, one has
\[
(h,f,g)=(-1)^{|f||g|}(h,g,f).
\]
\end{thm}
\begin{proof}
Use the Getzler identity.
\end{proof}

\begin{defn}
The \emph{derivation deviation} of $\de$ over the tribraces
$\{\cdot,\cdot,\cdot\}$ is defined by
\[
\dev_{\{\cdot,\cdot,\cdot\}}\de\,(h\t f\t g)
    \=\de\{h,f,g\}-\{h,f,\de g\}-(-1)^{|g|}\{h,\de f,g\}
      -(-1)^{|g|+|f|}\{\de h,f,g\}.
\]
\end{defn}

\begin{thm}[{\rm\cite{KP}}]
\label{tri-main}
In a pre-operad $C$, one has
\[
(-1)^{|g|}\dev_{\{\cdot,\cdot,\cdot\}}\de\,(h\t f\t g)=
    (h\bul f)\u g+(-1)^{|h|f}f\u(h\bul g)-h\bul(f\u g).
\]
\end{thm}

\begin{thm}
In a pre-operad $C$, one has
\[
(-1)^{|g|}\dev_{\{\cdot,\cdot,\cdot\}}\de\,(h\t f\t g)=
    [h,f]\u g+(-1)^{|h|f}f\u[h,g]-[h,f\u g].
\]
\end{thm}
\begin{proof}
Combine Theorem \ref{tri-main} with Theorem \ref{right der}.
\end{proof}

\section{Main theorem and Gerstenhaber's method}
\label{tetrabraces}

\noindent
In this section, we calculate the derivation deviation of $\de$ over
tetrabraces for an abstract pre-operad.
For this we use the Gerstehaber (auxiliary variables) method
\cite{Ger,KPS}.

\begin{defn}[ground tetrahedron]
The \emph{ground tetrahedron} $T$
(associated with $h\t f\t g\in C^h\t C^f\t C^g$) is defined as
a set of points in $\NN\x\NN\x\NN$ with $\NN$-valued coordinates
$i,j,k$, which satisfy inequalities
\[
\begin{cases}
   0\leq i\leq h-3,\\
i+f \leq j\leq h+f-3,\\
j+g \leq k\leq h+f+g-3.
\end{cases}
\]
In other words,
\[
T\=\{(i,j,k)\in\NN\x\NN\x\NN|\,0\leq i\leq j-f\leq k-f-g\leq h-3\}.
\]
\end{defn}

\begin{defn}[tetrabraces]
Let $h\t f\t g\t b\in C^h\t C^f\t C^g \t C^b$. Define the tetrabraces
$\{\cdot,\cdot,\cdot,\cdot\}$ by
\[
\{h,f,g,b\}=\sum_{(i,j,k)\in T}((h\o_{i}f)\o_{j}g)\o_{k}b\in C^{h+f+g+b-3},
\qquad|\{\cdot,\cdot,\cdot,\cdot\}|=0.
\]
\end{defn}

\begin{defn}
The \emph{derivation deviation} of $\de$ over the tetrabraces
$\{\cdot,\cdot,\cdot,\cdot\}$ is defined by
\begin{align*}
\dev_{\{\cdot,\cdot,\cdot,\cdot\}}\de\,(h\t f\t g\t b)
    &\=\de\{h,f,g,b\}-\{h,f,g,\de b\}-(-1)^{|b|}\{h,f,\de g,b\}\\
    &\quad-(-1)^{|b|+|g|}\{h,\de f,g,b\}
          -(-1)^{|b|+|g|+|f|}\{\de h,f,g,b\}.
\end{align*}
\end{defn}

\begin{Main Theorem}
In a pre-operad $C$, one has
\begin{align*}
(-1)^{|b|}\dev_{\{\cdot,\cdot,\cdot,\cdot\}}\de\,(h\t f\t g\t b)
    &=\{h,f,g\}\u b-\{h,f,g\u b\}\\
    &-(-1)^{|g|}\{h,f\u g,b\}+(-1)^{|h|f+|g|}f\u\{h,g,b\}.
\end{align*}
\end{Main Theorem}

To prove this theorem, we use auxiliary variables.

\begin{notations}[auxiliary variables]
\label{aux var}
In a pre-operad $C$, for $(i,j,k)\in T$ define the \emph{auxiliary variables}
\allowdisplaybreaks
\begin{align*}
\Ga_{i+1,j+1,k+1}
    &\=-(-1)^{|h|+|f|+|g|+|b|}\,\I\u(((h\o_i f)\o_j g)\o_k b)\\
    &\quad\,\,\,-(-1)^{|f|+|g|+|b|}\sum_{s=0}^{i-1}
         (((h\o_s\mu)\o_{i+1}f)\o_{j+1}g)\o_{k+1} b\\
    &\quad\,\,\,+(-1)^{|f|+|g|+|b|}((h\o_i(\I\u f))\o_{j+1}g)\o_{k+1}b, \\
\Ga'_{i+1,j+1,k+1}
       &\=(-1)^{|g|+|b|}((h\o_i(f\u\I))\o_{j+1} g)\o_{k+1}b\\
       &\quad\,\,\,-(-1)^{|f|+|g|+|b|}\sum_{s=i+1}^{j-f}
         (((h\o_s\mu)\o_i f)\o_{j+1}g)\o_{k+1}b\\
       &\quad\,\,\,+(-1)^{|g|+|b|}((h\o_i f)\o_{j}(\I\u g))\o_{k+1}b, \\
\Ga''_{i+1,j+1,k+1}
       &\=(-1)^{|b|}((h\o_i f)\o_j(g\u\I))\o_{k+1} b\\
       &\quad\,\,\,-(-1)^{|f|+|g|+|b|}\sum_{s=j-|f|+1}^{k-|f|-g}
         (((h\o_s\mu)\o_i f)\o_j g)\o_{k+1}b\\
       &\quad\,\,\,+(-1)^{|b|}((h\o_i f)\o_j g)\o_{k}(\I\u b), \\
\Ga'''_{i+1,j+1,k+1}
       &\=((h\o_i f)\o_j g)\o_k(b\u\I)\\
       &\quad\,\,\,-(-1)^{|f|+|g|+|b|}\sum_{s=k-|f|-|g|+1}^{|h|}
         (((h\o_s\mu)\o_if)\o_jg)\o_k b\\
       &\quad\,\,\,-(((h\o_i f)\o_j g)\o_k b)\u\I.
\end{align*}
\allowdisplaybreaks
\end{notations}

\begin{lem}
\label{first}
In a pre-operad $C$, for $(i,j,k)\in T$ one has
\begin{align*}
\de(((h\o_i f)\o_j g)\o_k b)
    &-((h\o_i f)\o_j g)\o_k\de b
         -(-1)^{|b|}((h\o_i f)\o_j\de g)\o_{k+1}b\\
    &-(-1)^{|b|+|g|}((h\o_i\de f)\o_{j+1}g)\o_{k+1} b\\
    &=\Ga_{i+1,j+1,k+1}+\Ga'_{i+1,j+1,k+1}+\Ga''_{i+1,j+1,k+1}
         +\Ga'''_{i+1,j+1,k+1}.
\end{align*}
\end{lem}
\begin{proof}
See Appendix A.
\end{proof}

\begin{notations}[truncated envelope of $T'$]
Now define a \emph{shifted} ground tetrahedron
\[
T':=\{(i,j,k)\in \NN\times\NN\times\NN|\,1\leq i\leq j-f\leq k-f-g\leq h-2\}.
\]
Its \emph{envelope} is the tetrahedron
\[
T'_{env}\=\{(i,j,k)\in\NN\times\NN\times\NN|\,
                 0\leq i\leq j-|f|\leq k-|f|-|g|\leq h+1\}.
\]
The \emph{boundary} of the envelope $T'_{env}$ is evidently
$\partial T'_{env}=T'_{env}\setminus T'$.
The \emph{truncated envelope} $T'_{\widetilde{env}}$ of $T'$ is defined by
\emph{removing} (six) edges of $T'_{env}$,
\begin{align*}
T'_{\widetilde{env}}
    &\=T'_{env}\setminus\{(i,j,k)\in\NN\x\NN\x\NN|\,(0\leq i\leq h+1,i+|f|,h+f+|g|);\\
    &\hskip5.37true cm(0,|f|,|f|+|g|\leq k\leq h+f+|g|);\\
    &\hskip5.37true cm(0\leq i\leq h+1,i+|f|,i+|f|+|g|);\\
    &\hskip5.37true cm(0,|f|\leq j\leq h+f, h+f+|g|);\\
    &\hskip5.37true cm(0\leq i\leq h+1,h+f,h+f+|g|);\\
    &\hskip5.37true cm(0,|f|\leq j\leq h+f,j+|g|)\}\\
    &\,\,=T'\sqcup\partial T'_{\widetilde{env}}.
\end{align*}
\end{notations}

\begin{lem}
\label{second}
In a pre-operad $C$, for $0\leq i\leq j-f\leq k-f-g\leq h-2$
one has
\[
(-1)^{|f|+|g|+|b|}((\de h\o_i f)\o_j g)\o_k b
    =\Ga_{ijk}+\Ga'_{i+1,jk}+\Ga''_{i+1,j+1,k}+\Ga'''_{i+1,j+1,k+1},
\]
by definition for $\Ga_{0jk}$, $\Ga'_{i,i+|f|,k}$,
$\Ga''_{ij,j+|g|}$, and $\Ga'''_{ij,|h|+f+g}$
(boundary values on $\partial T'_{\widetilde{env}}$).
\end{lem}
\begin{proof}
See Appendix B.
\end{proof}

\begin{Boundary Lemma}
In a pre-operad $C$, for $(i,j,k)\in \partial T'_{\widetilde{env}}$ one has
\allowdisplaybreaks
\begin{alignat*}{6}
&\Ga_{0jk}      &&=(-1)^{|g|+b+|h|f}f\u((h\o_{j-f}g)\o_{k-f}b),
    &&\quad f\leq j\leq k-g\leq |h|+|f|,\\
&\Ga'_{i,i+|f|,k}&&=(-1)^{|b|+|g|}   (h\o_{i-1}(f\u g))\o_{k} b,
    &&\quad 1\leq i\leq k-|f|-g\leq |h|,\\
&\Ga''_{ij,j+|g|}&&=(-1)^{|b|}   (h\o_{i-1}f)\o_{j-1}(g\u b),
    &&\quad 1\leq i\leq j-f\leq |h|,\\
&\Ga'''_{ij,|h|+f+g} &&=(-1)^{b}     ((h\o_{i-1}f)\o_{j-1}g)\u b,
    &&\quad 1\leq i\leq j-f\leq |h|,
\end{alignat*}
\allowdisplaybreaks
by definition for the six edge values
\begin{align*}
&\Ga_{0j,|h|+|f|+g}\,\quad(f\leq j\leq|h|+|f|),       &\quad
      &\Ga'_{1fk}\,\quad(f+g\leq k\leq|h|+|f|+|g|),\\
&\Ga'_{i,i+|f|,|h|+|f|+g}\,\quad(1\leq i\leq|h|),     &\quad
      &\Ga''_{1j,j+|g|}\,\quad(f+1\leq j\leq|h|+|f|),\\
&\Ga''_{i,i+f,i+f+|g|}\,\quad\quad(2\leq i\leq|h|-1), &\quad
      &\Ga''_{i,|h|+f,h+|f|+|g|}\,\quad(1\leq i\leq|h|).
\end{align*}
\end{Boundary Lemma}
\begin{proof}
See Appendix C.
\end{proof}

\subsection*{Proof of the Main Theorem (Gerstenhaber's method)}
First note that
\allowdisplaybreaks
\begin{align*}
\{h,f,\de g,b\}
 &=\sum_{i=0}^{|h|-2}\sum_{j=i+f}^{|h|+|f|-1}
    \sum_{k=j+\de g}^{|h|+|f|+|\de g|}((h\o_i f)\o_j\de g)\o_k b\\
 &=\sum_{i=0}^{|h|-2}\sum_{j=i+f}^{|h|+|f|-1}
    \sum_{k=j+g+1}^{|h|+|f|+|g|+1}((h\o_i f)\o_j\de g)\o_k b\\
 &=\sum_{i=0}^{|h|-2}\sum_{j=i+f}^{|h|+|f|-1}
    \sum_{k=j+g}^{|h|+|f|+|g|}((h\o_i f)\o_{j}\de g)\o_{k+1} b\\
 &=\sum_{(i,j,k)\in T}((h\o_i f)\o_{j}\de g)\o_{k+1} b
\end{align*}
\allowdisplaybreaks
and
\allowdisplaybreaks
\begin{align*}
\{h,\de f,g,b\}
 &=\sum_{i=0}^{|h|-2}\sum_{j=i+\de f}^{|h|+|\de f|-1}
    \sum_{k=j+g}^{|h|+|\de f|+|g|}((h\o_i \de f)\o_j g)\o_k b\\
 &=\sum_{i=0}^{|h|-2}\sum_{j=i+f+1}^{|h|+|f|}
    \sum_{k=j+g}^{|h|+|f|+|g|+1}((h\o_i \de f)\o_j g)\o_k b\\
 &=\sum_{i=0}^{|h|-2}\sum_{j=i+f}^{|h|+|f|-1}
    \sum_{k=j+g}^{|h|+|f|+|g|}((h\o_i \de f)\o_{j+1}g)\o_{k+1} b\\
 &=\sum_{(i,j,k)\in T}((h\o_i \de f)\o_{j+1}g)\o_{k+1} b.
\end{align*}
\allowdisplaybreaks
By using Lemma~\ref{first} we have
\allowdisplaybreaks
\allowdisplaybreaks
\begin{align*}
\de\{h,f,g,b\}&-\{h,f,g,\de b\}-(-1)^{|b|}\{h,f,\de g,b\}
    -(-1)^{|b|+|g|}\{h,\de f,g,b\}\\
    &=\sum_{(i,j,k)\in T}
     (\Ga_{i+1,j+1,k+1}+\Ga'_{i+1,j+1,k+1}+\Ga''_{i+1,j+1,k+1}
      +\Ga'''_{i+1,j+1,k+1})\\
&=\sum_{(i,j,k)\in T'}(\Ga_{ijk}+\Ga'_{ijk}+\Ga''_{ijk}+\Ga'''_{ijk}).
\end{align*}
\allowdisplaybreaks
Now use Lemma~\ref{second} to see that
\allowdisplaybreaks
\begin{align*}
(-1)^{|f|+|g|+|b|}&\{\de h,f,g,b\}
=(-1)^{|f|+|g|+|b|}\!\!\sum_{i=0}^{|\de h|-2}\sum_{j=i+f}^{|\de h|+|f|-1}
    \sum_{k=j+g}^{|\de h|+|f|+|g|}\!\!\!((\de h\o_i f)\o_j g)\o_k b\\
&= \sum_{i=0}^{|h|-1}\sum_{j=i+f}^{|h|+|f|}\sum_{k=j+g}^{h+|f|+|g|}
            (\Ga_{ijk}+\Ga'_{i+1,jk}+\Ga''_{i+1,j+1,k}
        +\Ga'''_{i+1,j+1,k+1})\\
&= \sum_{i=0}^{|h|-1}\sum_{j=i+f}^{|h|+|f|}\sum_{k=j+g}^{h+|f|+|g|}\Ga_{ijk}
  +\sum_{i=1}^{|h|}\sum_{j=i+|f|}^{|h|+|f|}\sum_{k=j+g}^{h+|f|+|g|}\Ga'_{ijk}\\
&\quad
  +\sum_{i=1}^{|h|}\sum_{j=i+f}^{h+|f|}\sum_{k=j+|g|}^{h+|f|+|g|}\Ga''_{ijk}
  + \sum_{i=1}^{|h|}\sum_{j=i+f}^{h+|f|}\sum_{k=j+g}^{h+f+|g|}\Ga'''_{ijk}\\
&= \sum_{i=1}^{|h|-1}\sum_{j=i+f}^{|h|+|f|}\sum_{k=j+g}^{h+|f|+|g|}\Ga_{ijk}
  +\sum_{j=f}^{|h|+|f|}\sum_{k=j+g}^{h+|f|+|g|}\Ga_{0jk}\\
&\quad
  +\sum_{i=1}^{|h|-1}\sum_{j=i+f}^{|h|+|f|}\sum_{k=j+g}^{h+|f|+|g|}\Ga'_{ijk}
  +\sum_{i=1}^{|h|}\sum_{k=i+|f|+g}^{h+|f|+|g|}\Ga'_{i,i+|f|,k}\\
& \quad
  +\sum_{i=1}^{|h|-1}\sum_{j=i+f}^{|h|+|f|}\sum_{k=j+g}^{h+|f|+|g|}\Ga''_{ijk}
  +\sum_{i=1}^{|h|}\sum_{j=i+f}^{|h|+f}\Ga''_{ij,j+|g|}\\
&\quad
  +\sum_{i=1}^{|h|-1}\sum_{j=i+f}^{|h|+|f|}\sum_{k=j+g}^{h+|f|+|g|}\Ga'''_{ijk}
  +\sum_{i=1}^{|h|}\sum_{j=i+f}^{|h|+f}\Ga'''_{ij,h+f+|g|}\\
&= \sum_{(i,j,k)\in T'}(\Ga_{ijk}+\Ga'_{ijk}+\Ga''_{ijk}+\Ga'''_{ijk})
\end{align*}
\vskip-0.2true cm
\hskip2.5cm
\framebox{\begin{minipage}{8true cm}
    \vspace{-\abovedisplayskip}
\begin{align*}
       +\sum_{j=f}^{|h|+|f|}\sum_{k=j+g}^{h+|f|+|g|}\Ga_{0jk}
    &+\sum_{i=1}^{|h|}\sum_{k=i+|f|+g}^{h+|f|+|g|}\Ga'_{i,i+|f|,k}\\
       +\sum_{i=1}^{|h|}\sum_{j=i+f}^{|h|+f}\Ga''_{ij,j+|g|}\,\,\,
    &+\sum_{i=1}^{|h|}\sum_{j=i+f}^{|h|+f}\Ga''_{ij,h+f+|g|}
\end{align*}
\end{minipage}}
\allowdisplaybreaks
\vskip11pt
\noindent
One can see that the resulting \boxed{\emph{boxed formula}} is a sum
over the \emph{boundary $\partial T'_{\widetilde{env}}$ of the truncated envelope}
$T'_{\widetilde{env}}$ of $T'$.
By cancelling the sums $\sum_{T'}$ and using the Boundary Lemma,
we finally obtain
\allowdisplaybreaks
\begin{align*}
\dev_{\{\cdot,\cdot,\cdot,\cdot\}} \de\,(h\t f\t g\t b)
&=(-1)^{|g|+|b|+|h|f}\sum_{j=f}^{|f|+|h|}\sum_{k=j+g}^{h+|f|+|g|}
            f\u((h\o_{j-f}g)\o_{k-f}b)\\
 &\quad-(-1)^{|g|+|b|} \sum_{i=1}^{|h|}\sum_{k=i+|f|+g}^{h+|f|+|g|}
            (h\o_{i-1}(f\u g))\o_{k}b\\
 &\quad-(-1)^{|b|} \sum_{i=1}^{|h|}\sum_{j=i+f}^{|h|+f}
            (h\o_{i-1} f)\o_{j-1}(g\u b)\\
 &\quad-(-1)^{b}  \sum_{i=1}^{|h|}\sum_{j=i+f}^{|h|+f}
            ((h\o_{i-1}f)\o_{j-1}g)\u b\\
=(-1)^{|b|}\Big[&(-1)^{|h|f+|g|}\sum_{j=0}^{|h|-1}\sum_{k=j+g}^{|h|+|g|}
        f\u((h\o_j g)\o_k b)\\
 &\quad-(-1)^{|g|}\sum_{i=0}^{|h|-1}\sum_{j=i+f+g}^{h+|f|+|g|}
        (h\o_i (f\u g))\o_{j}b\\
 -\sum_{i=0}^{|h|-1}\sum_{j=i+f}^{|h|+|f|}(h&\o_i f)\o_j (g\u b)
  +\sum_{i=0}^{|h|-1}\sum_{j=i+f}^{|h|+|f|}((h\o_i f)\o_j g)\u b\Big]\\
=(-1)^{|b|}\Big[(-1)^{|h|f+|g|}&f\u\{h,g,b\}-(-1)^{|g|}\{h,f\u g,b\}
-\{h,f,g\u b\}\\
&\quad+\{h,f,g\}\u b\Big],
\end{align*}
\allowdisplaybreaks
which is the required formula.
\qed

\section{Appendix A}

\subsection*{Proof of Lemma~\ref{first}}
First note that
\allowdisplaybreaks
\begin{align*}
-\de_\mu(((h&\o_i f)\o_j g)\o_k b)=
           (-1)^{|h|+|f|+|g|+|b|}\,\I\u(((h\o_i f)\o_j g)\o_k b)\\
           +&\sum_{s=0}^{|h|+|f|+|g|+|b|}(((h\o_i f)\o_j g)\o_k b)\o_s\mu
           +(((h\o_i f)\o_j g)\o_k b)\u\I.
\end{align*}
\allowdisplaybreaks
By using the composition relations
\[
(((h\o_i f)\o_j g)\o_k b)\o_s\mu=
\begin{cases}
 (-1)^{|b|}(((h\o_i f)\o_j g)\o_s\mu)\o_{k+1}b,
                              &\text{if $0\leq s\leq k-1$}  \\
 ((h\o_i f)\o_j g)\o_k(b\o_{s-k}\mu),
                              &\text{if $k\leq s\leq k+|b|$}\\
 (-1)^{|b|}(((h\o_i f)\o_j g)\o_{s-|b|}\mu)\o_k b,
                              &\text{if $k+b\leq s$}
\end{cases}
\]
cut the above sum $\sum_{s=0}^{|h|+|f|+|g|+|b|}$ in three pieces,
\allowdisplaybreaks
\begin{align*}
-\de(((&h\o_i f)\o_j g)\o_k b)
   =(-1)^{|h|+|f|+|g|+|b|}\,\I\u(((h\o_i f)\o_j g)\o_k b)\\
       &+(-1)^{|b|}\sum_{s=0}^{k-1}(((h\o_i f)\o_j g)\o_s\mu)\o_{k+1}b
          +\sum_{s=k}^{k+|b|}((h\o_i f)\o_j g)\o_k(b\o_{s-k}\mu)\\
       &+(-1)^{|b|}\sum_{s=k+b}^{|h|+|f|+|g|+|b|}
          (((h\o_i f)\o_j g)\o_{s-|b|}\mu)\o_k b
          +(((h\o_i f)\o_j g)\o_k b)\u\I.
\end{align*}
\allowdisplaybreaks
Next, cut the first sum $\sum_{s=0}^{k-1}$ in three pieces
by using the composition relations
\[
((h\o_i f)\o_j g)\o_s\mu=
\begin{cases}
 (-1)^{|g|}((h\o_i f)\o_s\mu)\o_{j+1}g,
                              &\text{if $0\leq s\leq j-1$}  \\
 (h\o_i f)\o_j(g\o_{s-j}\mu),
                              &\text{if $j\leq s\leq j+|g|$}\\
 (-1)^{|g|}((h\o_i f)\o_{s-|g|}\mu)\o_j g,
                              &\text{if $j+g\leq s\leq k-1$},
\end{cases}
\]
and obtain
\allowdisplaybreaks
\begin{align*}
-\de(((h\o_i f)\o_j g)\o_k b)
   =&\,\,(-1)^{|h|+|f|+|g|+|b|}\,\I\u(((h\o_i f)\o_j g)\o_k b)\\
    &\,\,+(-1)^{|b|+|g|}\sum_{s=0}^{j-1}(((h\o_i f)\o_s\mu)\o_{j+1}g)\o_{k+1}b\\
    &\,\,+(-1)^{|b|}\sum_{s=j}^{j+|g|}((h\o_i f)\o_j(g\o_{s-j}\mu))\o_{k+1}b\\
    &\,\,+(-1)^{|b|+|g|}\sum_{s=j+g}^{k-1}(((h\o_i f)\o_{s-|g|}\mu)\o_j g)\o_{k+1}b\\
    &\,\,+\sum_{s=k}^{k+|b|}((h\o_i f)\o_j g)\o_k(b\o_{s-k}\mu)\\
    &\,\,+(-1)^{|b|}\sum_{s=k+b}^{|h|+|f|+|g|+|b|}
        (((h\o_i f)\o_j g)\o_{s-|b|}\mu)\o_k b\\
    &\,\,+(((h\o_i f)\o_j g)\o_k b)\u\I.
\end{align*}
\allowdisplaybreaks
Next, cut the first sum $\sum_{s=0}^{j-1}$ in three pieces 
by using the composition relations
\[
(h\o_i f)\o_s\mu=
\begin{cases}
 (-1)^{|f|}(h\o_s\mu)\o_{i+1}f,    &\text{if $0\leq s\leq i-1$}  \\
 h\o_i(f\o_{s-i}\mu),              &\text{if $i\leq s\leq i+|f|$}\\
 (-1)^{|f|}(h\o_{s-|f|}\mu)\o_i f, &\text{if $i+f\leq s\leq j-1$},
\end{cases}
\]
and desuspend summation ranges,
\[
\sum_{s=n}^{n+|x|} x\o_{s-n}\mu
  =\sum_{s=0}^{|x|}x\o_s\mu
  =x\bul\mu
  =-(-1)^{|x|}\,\I\u x-\de_\mu x-x\u\I
\]
for the pairs $(n,x)=(k,b)$, $(j,g)$ and $(i,f)$,
\allowdisplaybreaks
\begin{gather*}
\sum_{s=i+f}^{j-1} h\o_{s-|f|}\mu
 =\sum_{s=i+1}^{j-f} h\o_{s}\mu,\\
\sum_{s=j+g}^{k-1}(h\o_i f)\o_{s-|g|}\mu
  =(-1)^{|f|}\!\!\!\sum_{s=j+g}^{k-1}\!\!\!(h\o_{s-|g|-|f|}\mu)\o_i f
  =(-1)^{|f|}\!\!\!\sum_{s=j-|f|+1}^{k-|f|-g}\!\!\!(h\o_{s}\mu)\o_i f,\\
\begin{align*}
  \sum_{s=k+b}^{|h|+|f|+|g|+|b|}((h\o_i f)&\o_j g)\o_{s-|b|}\mu
   =(-1)^{|g|}\sum_{s=k+b}^{|h|+|f|+|g|+|b|}((h\o_i f)\o_{s-|b|-|g|}\mu)\o_j g\\
  &=(-1)^{|f|+|g|}\sum_{s=k+b}^{|h|+|f|+|g|+|b|}((h\o_{s-|b|-|g|-|f|}\mu)\o_i f)\o_j g\\
  &=(-1)^{|f|+|g|}\sum_{s=k-|f|-|g|+1}^{|h|}((h\o_{s}\mu)\o_i f)\o_j g
\end{align*}\\
\end{gather*}
\allowdisplaybreaks
to obtain the required formula.
\qed

\section{Appendix B}
\subsection{Proof of Lemma~\ref{second}}
\label{ground}
First note that
\allowdisplaybreaks
\begin{align*}
\Ga_{ijk}+\Ga'_{i+1,jk}&+\Ga''_{i+1,j+1,k}+\Ga'''_{i+1,j+1,k+1}\\
   =&-(-1)^{|h|+|f|+|g|+|b|}\,\I\u(((h\o_{i-1}f)\o_{j-1}g)\o_{k-1}b)\\
    &-(-1)^{|f|+|g|+|b|}\sum_{s=0}^{i-2}(((h\o_s\mu)\o_{i}f)\o_{j}g)\o_{k}b\\
    &+(-1)^{|f|+|g|+|b|}((h\o_{i-1}(\I\u f))\o_{j}g)\o_k b \\
    &+(-1)^{|g|+|b|}((h\o_i(f\u\I))\o_{j}g)\o_k b\\
    &-(-1)^{|f|+|g|+|b|}\sum_{s=i+1}^{j-f-1}(((h\o_s\mu)\o_{i}f)\o_{j}g)\o_k b\\
    &+(-1)^{|g|+|b|}((h\o_i f)\o_{j-1}(\I\u g))\o_{k}b \\
    &+(-1)^{|b|}((h\o_i f)\o_j(g\u\I))\o_{k}b\\
    &-(-1)^{|f|+|g|+|b|}\sum_{s=j-|f|+1}^{k-f-g}(((h\o_s\mu)\o_{i}f)\o_{j}g)\o_k b\\
    &+(-1)^{|b|}((h\o_i f)\o_j g)\o_{k-1}(\I\u b)\\
    & +((h\o_i f)\o_j g)\o_k(b\u\I)\\
    &-(-1)^{|f|+|g|+|b|}\sum_{s=k-|f|-|g|+1}^{|h|}(((h\o_s\mu)\o_{i}f)\o_j g)\o_k b\\
    & -(((h\o_if)\o_jg)\o_k b)\u\I.
\end{align*}
\allowdisplaybreaks
We must compare it term by term with
\allowdisplaybreaks
\begin{align*}
-((\de h\o_i f)&\o_j g)\o_k b
 =\Bigg(\bigg(
      \Big((-1)^{|h|}\,\I\u h+\sum_{s=0}^{|h|}h\o_s\mu+h\u\I\Big)\o_i f
   \bigg)\o_j g\Bigg)\o_k b\\
 =&\,(-1)^{|h|}(((\I\u h)\o_i f)\o_j g)\o_k b
   +\sum_{s=0}^{i-2}(((h\o_s\mu)\o_i f)\o_j g)\o_k b \\
  & +(((h\o_{i-1}\mu)\o_i f)\o_j g)\o_k b
   +(((h\o_{i}\mu)\o_i f)\o_j g)\o_k b \\
  &+\sum_{s=i+1}^{j-f-1}(((h\o_s\mu)\o_i f)\o_j g)\o_k b
   +(((h\o_{j-f}\mu)\o_i f)\o_j g)\o_k b \\
  &+(((h\o_{j-|f|}\mu)\o_i f)\o_j g)\o_k g
   +\sum_{s=j-|f|+1}^{k-f-g}(((h\o_s\mu)\o_i f)\o_j g)\o_k b\\
  &+(((h\o_{k-|f|-g}\mu)\o_i f)\o_j g)\o_k b
   +(((h\o_{k-|f|-|g|}\mu)\o_i f)\o_j g)\o_k b\\
  & +\sum_{s=k-|f|-|g|+1}^{|h|}(((h\o_s\mu)\o_i f)\o_j g)\o_k b
   +(((h\u\I)\o_i f)\o_j g)\o_k b.
\end{align*}
\allowdisplaybreaks
Now, recall the sign $(-1)^{|f|+|g|+|b|}$ and use composition relations to
note the {\bf\emph{ground identities}}
\allowdisplaybreaks
\begin{align*}
(((\I\u h)\o_i f)\o_j g)\o_k b
      &=((\I\u(h\o_{i-1}f))\o_j g)\o_k b\\
      &=(\I\u((h\o_{i-1}f)\o_{j-1}g))\o_k b\\
      &=\I\u(((h\o_{i-1} f)\o_{j-1}g)\o_{k-1} b),\\
 (h\o_{i-1}\mu)\o_i f
      &=h\o_{i-1}(\mu\o_1 f)=-h\o_{i-1}(\I\u f), \\
 (h\o_{i}\mu)\o_i f
      &=h\o_{i}(\mu\o_0 f)=(-1)^fh\o_{i}(f\u\I),\\
 ((h\o_{j-f}\mu)\o_i f)\o_j g
      &= (-1)^{|f|}((h\o_{i}f)\o_{j-1}\mu)\o_j g\\
      &= (-1)^{|f|}(h\o_{i}f)\o_{j-1}(\mu\o_1 g)\\
      &= (-1)^{f}(h\o_{i}f)\o_{j-1}(\I\u g),\\
((h\o_{j-|f|}\mu)\o_i f)\o_j g
      &= (-1)^{|f|}((h\o_{i}f)\o_j\mu)\o_j g\\
      &= (-1)^{|f|}(h\o_{i}f) \o_j(\mu\o_0 g) \\
      &= (-1)^{|f|+g}(h\o_{i}f)\o_j(g\u\I),\\
 (((h\o_{k-|f|-g}\mu)\o_i f)\o_j g)\o_k b
      &= (-1)^{|f|}(((h\o_{i}f)\o_{k-g}\mu)\o_j g)\o_k b\\
      &= (-1)^{|f|+|g|}(((h\o_{i}f)\o_j g)\o_{k-1}\mu)\o_k b)\\
      &= (-1)^{|f|+|g|}((h\o_{i}f)\o_j g)\o_{k-1}(\mu\o_1 b)\\
      &= (-1)^{|f|+g}((h\o_{i}f)\o_j g)\o_{k-1}(\I\u b),\\
(((h\o_{k-|f|-|g|}\mu)\o_i f)\o_j g)\o_k b
      &= (-1)^{|f|}(((h\o_{i}f)\o_{k-|g|}\mu)\o_j g)\o_k b \\
      &= (-1)^{|f|+|g|}((h\o_{i}f) \o_j g)\o_k\mu)\o_k b \\
      &= (-1)^{|f|+|g|}((h\o_{i}f) \o_j g)\o_k(\mu\o_0 b) \\
      &= (-1)^{|f|+|g|+b}((h\o_{i}f)\o_j g)\o_k(b\u\I),\\
 (((h\u\I)\o_i f)\o_j g)\o_k b
      &=(-1)^{|f|}(((h\o_i f)\u\I)\o_j g)\o_k b\\
      &=(-1)^{|f|+|g|}(((h\o_i f)\o_j g)\u\I)\o_k b\\
      & =(-1)^{|f|+|g|+|b|}(((h\o_i f)\o_j g)\o_k b)\u\I,
\end{align*}
\allowdisplaybreaks
which prove the required formula.
\qed

\subsection{Proposition/recapitulation}
\label{recap}
In a pre-operad $C$, for
$(i,j,k)\in T'$ the auxiliary variables read
\begin{align*}
\Ga_{ijk}\=&-(-1)^{|h|+|f|+|g|+|b|}\,(((\I\u h)\o_i f)\o_j g)\o_k b\\
    &\, -(-1)^{|f|+|g|+|b|}\sum_{s=0}^{i-1}
        (((h\o_s\mu)\o_{i}f)\o_{j}g)\o_{k} b, \\
\Ga'_{ijk}\=&-(-1)^{|f|+|g|+|b|}\sum_{s=i-1}^{j-f}
        (((h\o_s\mu)\o_{i-1} f)\o_{j}g)\o_{k}b, \\
\Ga''_{ijk}\=&-(-1)^{|f|+|g|+|b|}\sum_{s=j-f}^{k-f-|g|}
        (((h\o_s\mu)\o_{i-1} f)\o_{j-1} g)\o_{k}b, \\
\Ga'''_{ijk}\=&-(-1)^{|f|+|g|+|b|}\sum_{s=k-f-|g|}^{|h|}
        (((h\o_s\mu)\o_{i-1}f)\o_{j-1}g)\o_{k-1} b\\
           &\,-(-1)^{|f|+|g|+|b|}(((h\u\I)\o_{i-1} f)\o_{j-1}g)\o_{k-1} b.
\end{align*}
\begin{proof}
Use the {\bf\emph{ground identities}} from the previous section \ref{ground}.
\end{proof}

\section{Appendix C}

\subsection*{Proof of the Boundary Lemma}
We calculate the boundary values of the auxiliary variables in a standard way,
by using Lemma \ref{second} and recapitulation formulae~\ref{recap}.
First prove that
\[
\Ga''_{ij,j+|g|}=(-1)^{|b|}(h\o_{i-1}f)\o_{j-1}(g\u b)
\]
for $2\leq i\leq|h|-2$, $i+f+1\leq j\leq|h|+|f|$.
Use Lemma \ref{second} and recapitulation formulae~\ref{recap} to note that
\begin{align*}
\Ga_{i-1,j-1,j+|g|}&+\Ga'_{i,j-1,j+|g|}+\Ga''_{ij,j+|g|}+\Ga'''_{ij,j+g}\\
     &=-(-1)^{|h|+|f|+|g|+|b|}\,(((\I\u h)\o_{i-1}f)\o_{j-1} g)\o_{j+|g|}b)\\
     &\quad-(-1)^{|f|+|g|+|b|}\sum_{s=0}^{i-2}
        (((h\o_s\mu)\o_{i-1}f)\o_{j-1}g)\o_{j+|g|} b\\
     &\quad-(-1)^{|f|+|g|+|b|}\sum_{s=i-1}^{j-f-1}
        (((h\o_s\mu)\o_{i-1}f)\o_{j-1}g)\o_{j+|g|}b
      +\Ga''_{ij,j+|g|}\\
     &\quad-(-1)^{|f|+|g|+|b|}\sum_{s=j-|f|}^{|h|}
        (((h\o_s\mu)\o_{i-1}f)\o_{j-1}g)\o_{j+|g|}b\\
     &\quad-(-1)^{|f|+|g|+|b|}(((h\u\I)\o_{i-1} f)\o_{j-1} g)\o_{j+|g|}b.
\end{align*}
We must compare it term by term with
\begin{align*}
-((&\de h\o_{i-1}f)\o_{j-1}g)\o_{j+|g|}b
           =(-1)^{|h|}(((\I\u h)\o_{i-1}f)\o_{j-1} g)\o_{j+|g|}b\\
     &+\sum_{s=0}^{|h|}(((h\o_s \mu)\o_{i-1} f)\o_{j-1} g)\o_{j+|g|}b
            +(((h\u\I)\o_{i-1}f)\o_{j-1} g)\o_{j+|g|}b\\
    =&\,\,(-1)^{|h|}(((\I\u h)\o_{i-1} f)\o_{j-1} g)\o_{j+|g|}b
            +\sum_{s=0}^{i-2}(((h\o_s \mu)\o_{i-1} f)\o_{j-1}g)\o_{j+|g|}b\\
     &+\sum_{s=i-1}^{j-f-1}(((h\o_s \mu)\o_{i-1}f)\o_{j-1}g)\o_{j+|g|}b
            +(((h\o_{j-f}\mu)\o_{i-1} f)\o_{j-1}g)\o_{j+|g|} b\\
     &+\sum_{s=j-|f|}^{|h|}(((h\o_s \mu)\o_{i-1}f)\o_{j-1}g)\o_{j+|g|}b
              +(((h\u\I)\o_{i-1} f)\o_{j-1}g)\o_{j+|g|}b.
\end{align*}
Now recall the sign $(-1)^{|f|+|g|+|b|}$ and use composition relations
to note that
\begin{align*}
(((h\o_{j-f}\mu)\o_{i-1} f)\o_{j-1}g)\o_{j+|g|}b
    &=(-1)^{|f|}(((h\o_{i-1} f)\o_{j-1}\mu)\o_{j-1}g)\o_{j+|g|}b\\
    &=(-1)^{|f|}((h\o_{i-1} f)\o_{j-1}(\mu\o_0 g))\o_{j+|g|}b\\
    &=(-1)^{|f|}(h\o_{i-1} f)\o_{j-1}((\mu\o_0 g)\o_g b)\\
    &=(-1)^{|f|+g}(h\o_{i-1}f)\o_{j-1}(g\u b),
\end{align*}
which lead one to the required formula for $\Ga''_{ij,j+|g|}$.
Also note that the last identities hold if
$1\leq i\leq|h|$, $i+f\leq j\leq|h|+f$.
The latter inequalities represent a \emph{projection}
on the $ij$-coordinate plane
of the following face of $T'_{\widetilde{env}}$:
\[
1\leq i\leq|h|,\qquad i+f\leq j\leq|h|+f,\qquad k=j+|g|.
\]
Therefore, in this case (projection) we have
\begin{align*}
\Ga''_{ij,j+|g|}
  &=(-1)^{|b|}(h\o_{i-1}f)\o_{j-1}(g\u b)\\
  &=-(-1)^{|f|+|g|+|b|}(((h\o_{j-f}\mu)\o_{i-1} f)\o_{j-1}g)\o_{j+|g|}b\\
  &=-(-1)^{|f|+|g|+|b|}\sum_{s=j-f}^{j-f}
        (((h\o_s\mu)\o_{i-1} f)\o_{j-1}g)\o_{j+|g|}b
\end{align*}
and the recapitulation formula
\[
\Ga''_{ijk}=-(-1)^{|f|+|g|+|b|}\sum_{s=j-f}^{k-f-|g|}
        (((h\o_s\mu)\o_{i-1} f)\o_{j-1} g)\o_{k}b
\]
holds in the tetrahedron
\[
1\leq i\leq|h|,\qquad i+f\leq j\leq|h|+f,\qquad j+|g|\leq k\leq|h|+f+|g|.
\]

Next prove that
\[
\Ga'_{i,i+|f|,k}=(-1)^{|b|+|g|}(h\o_{i-1}(f\u g))\o_{k}b
\]
for $2\leq i\leq|h|-1$, $i+|f|+g\leq k\leq|h|+|f|+|g|$.
Use Lemma \ref{second} and recapitulation formulae~\ref{recap} to note that
\allowdisplaybreaks
\begin{align*}
\Ga_{i-1,i+|f|,k}
 &+\Ga'_{i,i+|f|,k}+\Ga''_{i,i+f,k}+\Ga'''_{i,i+f,k+1}\\
=&-(-1)^{|h|+|f|+|g|+|b|}\,(((\I\u h)\o_{i-1} f)\o_{i+|f|} g)\o_{k} b)\\
 &-(-1)^{|f|+|g|+|b|}\sum_{s=0}^{i-2}
    (((h\o_s\mu)\o_{i-1}f)\o_{i+|f|}g)\o_{k} b
  +\Ga'_{i,i+|f|,k}\\
 &-(-1)^{|f|+|g|+|b|}\sum_{s=i}^{k-f-|g|}
    (((h\o_s\mu)\o_{i-1} f)\o_{i+|f|}g)\o_{k}b\\
 &-(-1)^{|f|+|g|+|b|}\sum_{s=k-|f|-|g|}^{|h|}
    (((h\o_s\mu)\o_{i-1} f)\o_{i+|f|} g)\o_{k}b\\
 &-(-1)^{|f|+|g|+|b|}
    (((h\u\I )\o_{i-1} f)\o_{i+|f|} g)\o_{k} b.
\end{align*}
\allowdisplaybreaks
We have to compare it term by term with
\begin{align*}
-((\de h&\o_{i-1}f)\o_{i+|f|} g)\o_{k} b=
            (-1)^{|h|}(((\I\u h)\o_{i-1} f)\o_{i+|f|} g)\o_{k} b\\
        &+\sum_{s=0}^{|h|}(((h\o_s \mu)\o_{i-1} f)\o_{i+|f|} g)\o_{k} b
            +(((h\u\I)\o_{i-1} f)\o_{i+|f|} g)\o_{k} b\\
       =&\,\,(-1)^{|h|}(((\I\u h)\o_{i-1} f)\o_{i+|f|} g)\o_{k} b
            +\sum_{s=0}^{i-2}(((h\o_s \mu)\o_{i-1} f)\o_{i+|f|} g)\o_{k} b\\
        &+(((h\o_{i-1} \mu)\o_{i-1} f)\o_{i+|f|} g)\o_{k} b
            +\sum_{s=i}^{k-f-|g|}(((h\o_s \mu)\o_{i-1} f)\o_{i+|f|} g)\o_{k} b\\
        &+\sum_{s=k-|f|-|g|}^{|h|}(((h\o_s \mu)\o_{i-1} f)\o_{i+|f|} g)\o_{k} b
            +(((h\u\I)\o_{i-1} f)\o_{i+|f|} g)\o_{k} b.
\end{align*}
Now recall the sign $(-1)^{|f|+|g|+|b|}$ and use composition relations
to note that
\begin{align*}
((h\o_{i-1}\mu)\o_{i-1}f)\o_{i+|f|}g
    &=(h\o_{i-1}(\mu\o_0 f))\o_{i-1+f}g\\
    &=h\o_{i-1}((\mu\o_{0}f)\o_{f}g)\\
    &=(-1)^{f}h\o_{i-1}(f\u g),
\end{align*}
which lead one to the required formula for $\Ga'_{i,i+|f|,k}$.
Also note that the last identities hold if
$1\leq i\leq|h|$, $i+|f|+g\leq k\leq|h|+|f|+g$.
The latter inequalities represent a \emph{projection}
on the $ik$-coordinate plane
of the following face of $T'_{\widetilde{env}}$:
\[
1\leq i\leq|h|,\qquad j=i+|f|,\qquad i+|f|+g\leq k\leq|h|+|f|+g.
\]
Therefore, in this case (projection) we have
\begin{align*}
\Ga'_{i,i+|f|,k}&=(-1)^{|b|+|g|}   (h\o_{i-1}(f\u g))\o_{k} b\\
&=-(-1)^{|f|+|g|+|b|}(((h\o_{i-1}\mu)\o_{i-1} f)\o_{i+|f|} g)\o_k b\\
&=-(-1)^{|f|+|g|+|b|}\sum_{s=i-1}^{i-1}(((h\o_s\mu)\o_{i-1} f)\o_{i+|f|} g)\o_k b
\end{align*}
and the recapitulation formula
\[
\Ga'_{ijk}=-(-1)^{|f|+|g|+|b|}\sum_{s=i-1}^{j-f}
        (((h\o_s\mu)\o_{i-1} f)\o_{j}g)\o_{k}b
\]
holds in the tetrahedron
\[
1\leq i\leq|h|,\qquad i+|f|\leq j\leq|h|+|f|,\qquad j+g\leq k\leq|h|+|f|+g.
\]

Next prove that
\[
\Ga_{0jk}=(-1)^{|g|+b+|h|f}f\u((h\o_{j-f}g)\o_{k-f}b)
\]
for $f\leq j\leq|h|+|f|-1$, $j+g\leq k\leq|h|+|f|+|g|$.
Use Lemma \ref{second} and recapitulation formulae~\ref{recap} to note that
\begin{align*}
\Ga_{0jk}&+\Ga'_{1j k}+\Ga''_{1,j+1,k}+\Ga'''_{1,j+1,k+1}=\Ga_{0jk}\\
     &-(-1)^{|f|+|g|+|b|}\sum_{s=0}^{j-f}
        (((h\o_s\mu)\o_{0}f)\o_{j}g)\o_{k} b\\
     &-(-1)^{|f|+|g|+|b|}\sum_{s=j-|f|}^{k-f-|g|}
        (((h\o_s\mu)\o_0 f)\o_{j}g)\o_{k}b\\
     &-(-1)^{|f|+|g|+|b|}\sum_{s=k-|f|-|g|}^{|h|}
        (((h\o_s\mu)\o_0 f)\o_{j} g)\o_{k}b\\
     &-(-1)^{|f|+|g|+|b|}
        (((h\u\I)\o_0 f)\o_{j} g)\o_{k} b.
\end{align*}
We must compare it term by term with
\begin{align*}
-((\de h\o_0 f)&\o_{j} g)\o_{k} b
                  =(-1)^{|h|}(((\I\u h)\o_0 f)\o_{j} g)\o_{k}b\\
               &\quad+\sum_{s=0}^{|h|}(((h\o_s \mu)\o_0 f)\o_{j} g)\o_{k}b
                  +(((h\u\I)\o_0 f)\o_{j}g)\o_{k} b\\
               &=(-1)^{|h|}(((\I\u h)\o_0 f)\o_{j} g)\o_{k}b
                  +\sum_{s=0}^{j-f}(((h\o_s\mu)\o_0 f)\o_{j} g)\o_{k}b\\
               &\quad+\sum_{s=j-|f|}^{k-f-|g|}(((h\o_s\mu)\o_0 f)\o_{j} g)\o_{k}b\\
               &\quad+\sum_{s=k-|f|-|g|}^{|h|}(((h\o_s\mu)\o_0 f)\o_{j} g)\o_{k}b
                  +(((h\u\I)\o_0 f)\o_{j}g)\o_{k}b.
\end{align*}
Now recall the sign $(-1)^{|f|+|g|+|b|}$ and use
composition relations to note that
\begin{align*}
    (-1)^{|h|}(((\I\u h)\o_0 f)\o_{j} g)\o_{k} b
    &=-(-1)^{|h|}(((\mu\o_1 h)\o_0 f)\o_{j} g)\o_{k} b\\
    &=-(-1)^{|h|+|h||f|}(((\mu\o_0 f)\o_{f} h)\o_j g)\o_k b\\
    &=-(-1)^{f|h|}((\mu\o_{0}f)\o_{f}(h\o_{j-f}g))\o_k b\\
    &=-(-1)^{f|h|}(\mu\o_0 f)\o_f((h\o_{j-f}g)\o_{k-f}b)\\
    &=(-1)^{f|h|+|f|} f\u ((h\o_{j-f} g)\o_{k-f} b),
\end{align*}
which leads one to the required formula for $\Ga_{0jk}$.
Also note that the last identities hold if
$f\leq j\leq|h|+|f|$, $j+g\leq k\leq|h|+|f|+g$.
The latter inequalities represent a \emph{projection}
on the $jk$-coordinate plane
of the following face of $T'_{\widetilde{env}}$:
\[
i=0,\qquad f\leq j\leq|h|+|f|,\qquad j+g\leq k\leq|h|+|f|+g.
\]
Therefore, in this case (projection) we have
\begin{align*}
\Ga_{0jk}&=(-1)^{|g|+b+|h|f}f\u ((h\o_{j-f}g)\o_{k-f}b)\\
         &=-(-1)^{|h|+|f|+|g|+|b|}(((\I\u h)\o_0 f)\o_{j} g)\o_{k}b.
\end{align*}

At last prove that
\[
\Ga'''_{ij,|h|+f+g} =(-1)^{b}((h\o_{i-1}f)\o_{j-1}g)\u b
\]
for $1\leq i \leq |h|$, $i+f\leq j\leq |h|+f$.
Use Lemma \ref{second} and recapitulation formulae~\ref{recap}.
If $i=1$, then
\begin{align*}
\Ga_{0,j-1,|h|+|f|+g}
       &+\Ga'_{1,j-1,|h|+|f|+g}
         +\Ga''_{1j,|h|+|f|+g}+\Ga'''_{1j,|h|+f+g}\\
   =-(-&1)^{|h|+|f|+|g|+|b|}
        (((\I\u h)\o_0 f)\o_{j-1} g)\o_{|h|+|f|+g} b\\
    -(-&1)^{|f|+|g|+|b|}\sum_{s=0}^{j-f-1}
        (((h\o_s\mu)\o_{0} f)\o_{j-1}g)\o_{|h|+|f|+g}b\\
    -(-&1)^{|f|+|g|+|b|}\sum_{s=j-f}^{|h|}
        (((h\o_s\mu)\o_{0}f)\o_{j-1}g)\o_{|h|+|f|+g}b
        +\Ga'''_{1j,|h|+f+g}\\
   =-(-&1)^{|h|+|f|+|g|+|b|}
        (((\I\u h)\o_0 f)\o_{j-1}g)\o_{|h|+|f|+g}b\\
    -(-&1)^{|f|+|g|+|b|}\sum_{s=0}^{|h|}
        (((h\o_s \mu)\o_{0}f)\o_{j-1}g)\o_{|h|+|f|+g}b
        +\Ga'''_{1j,|h|+f+g}.
\end{align*}
If $2\leq i\leq |h|$, then
\allowdisplaybreaks
\begin{align*}
\Ga_{i-1,j-1,|h|+|f|+g}
      &+\Ga'_{i,j-1,|h|+|f|+g}
          +\Ga''_{ij,|h|+|f|+g}+\Ga'''_{ij,|h|+f+g}\\
  =-(-&1)^{|h|+|f|+|g|+|b|}
        (((\I\u h)\o_{i-1}f)\o_{j-1}g)\o_{|h|+|f|+g}b\\
   -(-&1)^{|f|+|g|+|b|}\sum_{s=0}^{i-2}
        (((h\o_s \mu)\o_{i-1}f)\o_{j-1}g)\o_{|h|+|f|+g}b\\
   -(-&1)^{|f|+|g|+|b|}\sum_{s=i-1}^{j-f-1}
        (((h\o_s\mu)\o_{i-1} f)\o_{j-1}g)\o_{|h|+|f|+g}b\\
   -(-&1)^{|f|+|g|+|b|}\sum_{s=j-f}^{|h|}
        (((h\o_s\mu)\o_{i-1}f)\o_{j-1}g)\o_{|h|+|f|+g}b
    +\Ga'''_{ij,|h|+f+g}\\
  =-(-&1)^{|h|+|f|+|g|+|b|}
        (((\I\u h)\o_{i-1}f)\o_{j-1}g)\o_{|h|+|f|+g}b\\
   -(-&1)^{|f|+|g|+|b|}\sum_{s=0}^{|h|}
        (((h\o_s \mu)\o_{i-1}f)\o_{j-1}g)\o_{|h|+|f|+g}b
   +\Ga'''_{ij,|h|+f+g}.\\
\end{align*}
\allowdisplaybreaks
So the both cases ($i=1$ and $2\leq i\leq|h|$) can be
obtained from unique formula for $1\leq i\leq|h|$.
We must compare it term by term with
\begin{align*}
((\de&h\o_{i-1}f)\o_{j-1} g)\o_{|h|+|f|+g} b=
           -(-1)^{|h|}(((\I\u h)\o_{i-1} f)\o_{j-1} g)\o_{|h|+|f|+g} b\\
     &-\sum_{s=0}^{|h|}(((h\o_s \mu)\o_{i-1} f)\o_{j-1} g)\o_{|h|+|f|+g}b
           -(((h\u\I)\o_{i-1} f)\o_{j-1} g)\o_{|h|+|f|+g}b.
\end{align*}
Now recall the sign $(-1)^{|f|+|g|+|b|}$ and use composition relations
to note that
\begin{align*}
-(((h\u\I)\o_{i-1}f)\o_{j-1} g)\o_{|h|+|f|+g}b
&=(-1)^{|h|}(((\mu\o_0h)\o_{i-1} f)\o_{j-1} g)\o_{h+|f|+|g|}b\\
&=(-1)^{|h|}(\mu\o_0((h\o_{i-1} f)\o_{j-1} g))\o_{h+|f|+|g|}b\\
&=-(-1)^{|f|+|g|}((h\o_{i-1}f)\o_{j-1} g)\u b,
\end{align*}
which lead one to the required formula for $\Ga'''_{ij,|h|+f+g}$.
\qed

\vskip1cm
\noindent
Department of Mathematics,
Tallinn Technical University\\
Ehitajate tee 5, 19086 Tallinn, Estonia\\
e-mails: liivi.kluge.ttu@mail.ee and epaal@edu.ttu.ee

\end{document}